\documentclass[12pt,reqno]{amsart}

\usepackage{amssymb}
\usepackage{amsfonts}
\usepackage{amsthm}

\theoremstyle{plain} 
\newtheorem{theorem}{Theorem}

\newtheorem{lemma}[theorem]{Lemma}
\theoremstyle{definition} 
\newtheorem{definition}[theorem]{Definition}
\theoremstyle{remark}

\newcommand{\T}{\ensuremath{\mathbb{T}}}

\DeclareMathOperator{\czero}{C_{rd}}
\DeclareMathOperator{\crdone}{C^{\Delta}_{rd}}
\DeclareMathOperator{\crdtwo}{C^{\Delta^2}_{rd}}

\numberwithin{equation}{section}
\numberwithin{theorem}{section}

\small\normalsize

\begin{document}
\author[anderson]{Douglas R. Anderson}
\title[hyers--ulam stability]{Hyers--Ulam stability of second-order linear dynamic equations on time scales}
\address{Department of Mathematics and Computer Science, Concordia College, Moorhead, MN 56562 USA}
\email{andersod@cord.edu}
\urladdr{http://www.cord.edu/faculty/andersod/bib.html}


\begin{abstract}
We establish the stability of second-order linear dynamic equations on time scales in the sense of Hyers and Ulam. To wit, if an approximate solution of the second-order linear equation exists, then there exists an exact solution to the dynamic equation that is close to the approximate one.
\end{abstract}

\maketitle

\thispagestyle{empty}


\section{introduction}
In 1940, Ulam posed the following problem concerning the stability of functional equations: give conditions in order for a linear mapping near an approximately linear mapping to exist. The problem for the case of approximately additive mappings was solved by Hyers who proved that the
Cauchy equation is stable in Banach spaces, and the result of Hyers was generalized by Rassias.

Throughout this work we assume the reader has a working knowledge of time scales.


\begin{definition}
For real constants $\alpha$ and $\beta$, consider the second-order linear dynamic equation
\begin{equation}\label{2ndcc}
 x^{\Delta\Delta}(t)+\alpha x^{\Delta}(t)+\beta x(t)=0, \quad t\in[a,b]_\T.
\end{equation}
If whenever $y\in\crdtwo[a,b]{_\T}$ satisfies
$$ \left| y^{\Delta\Delta}+\alpha y^{\Delta}+\beta y \right|\le \varepsilon $$
on $[a,b]_\T$, there exists a solution $u\in\crdtwo[a,b]{_\T}$ of \eqref{2ndcc} such that $|y-u|\le K\varepsilon$ on $[a,b]_\T$ for some constant $K>0$, then \eqref{2ndcc} has Hyers-Ulam stability $[a,b]_\T$.
\end{definition}


\begin{theorem}[Constant Coefficients]
If the characteristic equation $\lambda^2+\alpha\lambda+\beta=0$ has two distinct positive roots, then \eqref{2ndcc} has Hyers-Ulam stability on $[a,b]_\T$.
\end{theorem}

\begin{proof}
Let $\varepsilon>0$ be given, and let $y\in\crdtwo[a,b]_\T$ such that $\left| y^{\Delta\Delta}+\alpha y^{\Delta}+\beta y \right|\le \varepsilon$ on $[a,b]_\T$. Let $\lambda_1,\lambda_2$ be the distinct positive roots of $\lambda^2+\alpha\lambda+\beta=0$. On $[a,b]_\T$ define
$$ g:=y^\Delta-\lambda_1 y; $$
then $g^\Delta=y^{\Delta\Delta}-\lambda_1 y^\Delta$, so that
\begin{eqnarray*}
 |g^\Delta-\lambda_2 g| = |y^{\Delta\Delta}-\lambda_1 y^\Delta-\lambda_2 y^\Delta+\lambda_1\lambda_2 y| = |y^{\Delta\Delta}+\alpha y^{\Delta}+\beta y|\le \varepsilon
\end{eqnarray*}
on $[a,b]_\T$. Thus $-\varepsilon \le g^\Delta-\lambda_2 g\le \varepsilon$, or rewritten,
$$ \frac{-\varepsilon}{1+\mu \lambda_2} \le g^\Delta+(\ominus\lambda_2) g^\sigma \le \frac{\varepsilon}{1+\mu\lambda_2}. $$
For the case $0<\lambda_2\le 1$ there exists $M>0$ such that $M\lambda_2>1$, so without loss of generality we can assume that $\lambda_2>1$. Then
$$ \varepsilon(\ominus\lambda_2) \le g^\Delta+(\ominus\lambda_2) g^\sigma \le -\varepsilon(\ominus\lambda_2). $$
Multiply by $e_{\ominus\lambda_2}(\cdot,a)$ to see that
$$ \varepsilon\left(e_{\ominus\lambda_2}(\cdot,a)\right)^\Delta(t)\le \left(ge_{\ominus\lambda_2}(\cdot,a)\right)^\Delta(t) \le -\varepsilon\left(e_{\ominus\lambda_2}(\cdot,a)\right)^\Delta(t), $$
so that delta integrating from $t$ to $b$ yields
$$ -\varepsilon\left(e_{\ominus\lambda_2}(t,a)-e_{\ominus\lambda_2}(b,a)\right) \le g(b)e_{\ominus\lambda_2}(b,a)-g(t)e_{\ominus\lambda_2}(t,a) \le \varepsilon\left(e_{\ominus\lambda_2}(t,a)-e_{\ominus\lambda_2}(b,a)\right). $$
Then
$$ -\varepsilon e_{\ominus\lambda_2}(t,a) \le \left(g(b)-\varepsilon\right)e_{\ominus\lambda_2}(b,a)-g(t)e_{\ominus\lambda_2}(t,a) \le \varepsilon e_{\ominus\lambda_2}(t,a)-2e_{\ominus\lambda_2}(b,a), $$
whence
$$ -\varepsilon e_{\ominus\lambda_2}(t,a) \le \left(g(b)-\varepsilon\right)e_{\ominus\lambda_2}(b,a)-g(t)e_{\ominus\lambda_2}(t,a) \le \varepsilon e_{\ominus\lambda_2}(t,a). $$
Multiplying the above inequality by $e_{\lambda_2}(t,a)$ results in
$$ -\varepsilon \le \left(g(b)-\varepsilon\right) e_{\ominus\lambda_2}(b,t)-g(t) \le \varepsilon. $$
If we let 
$$ z(t):=\left(g(b)-\varepsilon\right) e_{\ominus\lambda_2}(b,t)=\left(g(b)-\varepsilon\right) e_{\lambda_2}(t,b), \quad t\in[a,b]_\T, $$
then clearly $z^\Delta(t)=\lambda_2 z(t)$ and $|g(t)-z(t)|\le\varepsilon$ for $t\in[a,b]_\T$.
Since $g(t)=y^\Delta(t)-\lambda_1 y(t)$ for $t\in[a,b]_\T$, we have
$$ -\varepsilon \le y^\Delta(t)-\lambda_1 y(t)-z(t)\le \varepsilon. $$
By an argument similar to the one given above, we can show that there exists
$$ u(t):=(y(b)-\varepsilon)e_{\lambda_1}(t,b)-e_{\lambda_1}(t,a)\int_t^b\frac{z(s)}{1+\mu(s)\lambda_1}e_{\ominus\lambda_1}(s,a)\Delta s $$
such that $|y(t)-u(t)|\le\varepsilon$ for $t\in[a,b]_\T$ and $u\in\crdtwo[a,b]_\T$ satisfies $u^\Delta-\lambda_1 u-z=0$ on $[a,b]_\T$. Consequently
$$ u^{\Delta\Delta}(t)-(\lambda_1+\lambda_2)u^\Delta(t)+\lambda_1\lambda_2 u(t)=0, \quad t\in[a,b]_\T $$
that is
$$ u^{\Delta\Delta}+\alpha u^{\Delta}+\beta u=0 $$
on $[a,b]_\T$, completing the proof.
\end{proof}

For the next result consider the inhomogeneous second-order linear dynamic equation
\begin{equation}\label{2ndicc}
 x^{\Delta\Delta}(t)+\alpha x^{\Delta}(t)+\beta x(t) = f(t), \quad t\in[a,b]_\T.
\end{equation}


\begin{theorem}[Inhomogeneous with Constant Coefficients]
Assume the characteristic equation $\lambda^2+\alpha\lambda+\beta=0$ has two distinct positive roots. For every $\varepsilon>0$, $f\in\czero [a,b]_\T$, and $y\in\crdtwo [a,b]_\T$, if 
\begin{equation}
 |y^{\Delta\Delta}+\alpha y^{\Delta}+\beta y - f| \le \varepsilon 
\end{equation}
on $[a,b]_\T$, then there exists a solution $u\in\crdtwo[a,b]{_\T}$ of \eqref{2ndicc} such that $|y-u|\le K\varepsilon$ on $[a,b]_\T$ for some constant $K>0$, that is to say \eqref{2ndicc} has Hyers-Ulam stability on $[a,b]_\T$.
\end{theorem}

The next theorem considers the inhomogeneous second-order linear dynamic equation with variable coefficients
\begin{equation}\label{2ndivc}
 x^{\Delta\Delta}(t) + p(t)x^{\Delta}(t) + q(t)x(t) = f(t), \quad t\in[a,b]_\T.
\end{equation}
First we will need the following lemma.


\begin{lemma}\label{lemma14}
Let $d,f\in\czero[a,b]_\T$ such that $1+\mu(t)d(t)\neq 0$ for all $t\in[a,b]_\T$ and 
$$ \sup_{t\in[a,b]_\T} \left|e_d(t,a)\right| \int_a^t \left|e_d(a,\sigma(s))\right|\Delta s<\infty. $$
Let $x\in\crdone[a,b]_\T$. Then the first-order dynamic equation
\begin{equation}\label{1stvc}
 x^\Delta(t)-d(t)x(t)-f(t)=0, \quad t\in[a,b]_\T
\end{equation}
has Hyers-Ulam stability, that is whenever $g\in\crdone[a,b]{_\T}$ satisfies
$$ \left|g^\Delta(t)-d(t)g(t)-f(t)\right|\le \varepsilon, \quad t\in[a,b]_\T $$
there exists a solution $w\in\crdone[a,b]{_\T}$ of \eqref{1stvc} such that $|g-w|\le L\varepsilon$ on $[a,b]_\T$ for some constant $L>0$.
\end{lemma}

\begin{proof}
Given $\varepsilon>0$, suppose there exists $g\in\crdone[a,b]{_\T}$ that satisfies
$$ \left|g^\Delta(t)-d(t)g(t)-f(t)\right|\le \varepsilon, \quad t\in[a,b]_\T. $$
Set 
$$ \ell:=g^\Delta-dg-f; $$
by \cite[Theorem 2.77]{bp1} we have that $g$ is given by
$$ g(t)=e_d(t,a)g(a)+\int_a^t e_d(t,\sigma(s))\left(f(s)+\ell(s)\right)\Delta s. $$
Let $w$ be the unique solution of the initial value problem
$$ w^\Delta-dw-f=0, \quad w(a)=g(a). $$
Then
$$ w(t)=e_d(t,a)g(a)+\int_a^t e_d(t,\sigma(s))f(s)\Delta s, $$
and 
\begin{eqnarray*}
 |g(t)-w(t)| &=& \left|\int_a^t e_d(t,\sigma(s))\ell(s)\Delta s\right| \\
 &\le& \left|e_d(t,a)\int_a^t e_d(a,\sigma(s))\ell(s)\Delta s\right| \\
 &\le& \varepsilon \sup_{t\in[a,b]_\T} \left|e_d(t,a)\right| \int_a^t \left|e_d(a,\sigma(s))\right|\Delta s \\
 &\le& L\varepsilon
\end{eqnarray*} 
for all $t\in[a,b]_\T$, where $L:=\sup_{t\in[a,b]_\T} \left|e_d(t,a)\right| \int_a^t \left|e_d(a,\sigma(s))\right|\Delta s$ is a constant independent of $g$ and $\varepsilon$. Since $w$ solves \eqref{1stvc} by construction, the proof in complete.
\end{proof}


\begin{theorem}[Inhomogeneous with Variable Coefficients]
Let $p,q,f\in\czero[a,b]_\T$ and consider \eqref{2ndivc}. Assume the related dynamic Riccati equation
$$ z^\Delta(t)+p(t)z(t)-z(t)z^\sigma(t)=q(t), \quad t\in[a,b]_\T $$
has a particular solution $z$ with both $1+\mu(t)(z^\sigma(t)-p(t))\neq 0$ and $1-\mu(t)z(t)\neq 0$ for all $t\in[a,b]_\T$. Furthermore assume that
\begin{equation}\label{ivcz1}
 \sup_{t\in[a,b]_\T} \left|e_{(z^\sigma-p)}(t,a)\right| \int_a^t \left|e_{(z^\sigma-p)}(a,\sigma(s))\right|\Delta s<\infty 
\end{equation}
and
\begin{equation}\label{ivcz2} 
 \sup_{t\in[a,b]_\T} \left|e_{-z}(t,a)\right| \int_a^t \left|e_{-z}(a,\sigma(s))\right|\Delta s<\infty. 
\end{equation}
Then \eqref{2ndivc} has Hyers-Ulam stability on $[a,b]_\T$.
\end{theorem}

\begin{proof}
We need to show that if there exists a $y\in\crdtwo[a,b]_\T$ that satisfies 
\begin{equation}
 |y^{\Delta\Delta}(t)+p(t)y^{\Delta}(t)+q(t)y(t) - f(t)| \le \varepsilon 
\end{equation}
for $t\in[a,b]_\T$, and the dynamic Riccati equation
$$ z^\Delta(t)+p(t)z(t)-z(t)z^\sigma(t)=q(t) $$
has a particular solution $z$ with both $1+\mu(t)(z^\sigma(t)-p(t))\neq 0$ and $1-\mu(t)z(t)\neq 0$ for $t\in[a,b]_\T$ such that \eqref{ivcz1} and \eqref{ivcz2} hold, then there exists a solution $u\in\crdtwo[a,b]{_\T}$ of \eqref{2ndivc} such that $|y-u|\le K\varepsilon$ on $[a,b]_\T$ for some constant $K>0$.

Let $\varepsilon>0$ be given, and let $y\in\crdtwo[a,b]_\T$ such that $\left|y^{\Delta\Delta} + py^{\Delta} + qy - f\right|\le \varepsilon$ on $[a,b]_\T$. Assume $z$ is a particular solution of the Riccati equation $z^\Delta+pz-zz^\sigma=q$ on $[a,b]_\T$, and set
$$ g:=y^\Delta+zy, \qquad d:=z^\sigma-p. $$
Then $g^\Delta = y^{\Delta\Delta} + z^\sigma y^\Delta + z^\Delta y$, so that
\begin{eqnarray*}
 |g^\Delta-dg-f| = |y^{\Delta\Delta} + z^\sigma y^\Delta + z^\Delta y - (z^\sigma-p)(y^\Delta+zy) - f| = |y^{\Delta\Delta}+py^{\Delta}+qy-f|\le \varepsilon
\end{eqnarray*}
on $[a,b]_\T$. As all of the hypotheses of Lemma \ref{lemma14} hold, equation \eqref{1stvc} has Hyers-Ulam stability, and there exists a solution $w\in\crdone[a,b]{_\T}$ of
\begin{equation}\label{bp230}
 w^\Delta(t) - d(t)w(t)-f(t)=0, \quad t\in[a,b]_\T
\end{equation}
where $w$ is given by
$$ w(t)=e_d(t,a)g(a)+\int_a^t e_d(t,\sigma(s))f(s)\Delta s, $$
and there exists an $L>0$ such that
$$ |g(t)-w(t)|\le L\varepsilon, \quad t\in[a,b]_\T. $$
Since $g=y^\Delta+zy$, we have that
$$ |y^\Delta(t)+z(t)y(t)-w(t)|\le L\varepsilon, \quad t\in[a,b]_\T. $$
Again apply Lemma \ref{lemma14} to see that there exists a solution $u\in\crdone[a,b]{_\T}$ of
\begin{equation}\label{u230}
 u^\Delta(t) + z(t)u(t)-w(t)=0, \quad t\in[a,b]_\T
\end{equation}
given by
$$ u(t)=e_{-z}(t,a)y(a)+\int_a^t e_{-z}(t,\sigma(s))w(s)\Delta s, $$
and there exists an $K>0$ such that
$$ |y(t)-u(t)|\le KL\varepsilon, \quad t\in[a,b]_\T. $$
Moreover,
\begin{eqnarray*}
 u^{\Delta\Delta} + pu^{\Delta} + qu - f 
 &=& w^\Delta-z^\sigma u^\Delta-z^\Delta u+ pu^{\Delta} + qu - f \\
 &=& (dw+f) - (d+p)u^\Delta + (q-z^\Delta)u + pu^{\Delta} - f \\
 &=& d(w-u^\Delta-zu) \\
 &=& 0
\end{eqnarray*}
on $[a,b]_\T$, so that $u$ is a solution of \eqref{2ndivc}, and actually $u\in\crdtwo[a,b]{_\T}$.
\end{proof}


\end{document}